\documentclass[english,11pt]{article}
\usepackage{amsfonts}
\usepackage[latin1]{inputenc}
\usepackage{amsmath,amsthm,amssymb}
\usepackage{graphicx}
\usepackage{color}

\def\neweq#1{\begin{equation}\label{#1}}
\def\endeq{\end{equation}}

\newtheorem{theorem}{Theorem}[section]
\newtheorem{lemma}{Lemma}[section]
\newtheorem{definition}{Definition}[section]

\newtheorem{remark}{Remark}[section]
\input{tcilatex}
\begin{document}

\title{\textbf{Existence and Non-Existence Criteria for Solutions of a Schrö%
dinger Quasilinear System type}}
\author{Dragos-Patru Covei\break \\
{\small \ Constantin Brancusi University of \ Tg-Jiu, Bld. Republicii 1,
210152, Romania}\\
{\small E-mail: c\texttt{ovdra@yahoo.com}}}
\date{}
\maketitle

\begin{abstract}
In this paper we analyze the existence of large positive radial solutions to
some quasilinear elliptic systems. Also, a non-radially symmetric solution
is obtained by using a lower and upper solution method. The equations are
coupled by functions which are increasing with respect to all the variables. 
\\[3pt]
\end{abstract}

\baselineskip16pt \renewcommand{\theequation}{\arabic{section}.%
\arabic{equation}} \catcode`@=11 \@addtoreset{equation}{section} \catcode%
`@=12



\section{Introduction}

In the present work we prove some existence and non-existence theorems for
the solutions of a quasilinear elliptic pde's system of \textit{Schrödinger }%
type. Let $m\in \left\{ 1,2,3,...\right\} $ and $r:=\left\vert x\right\vert $
the Euclidean norm of $x\in \mathbb{R}^{N}$. We begin now the precise
statements of our existence theorems by assuming that the class of functions 
$a_{j}$, $f_{j}$ ($j=1,...,m$)\ \ satisfy:

(A1)\quad $a_{j}:\mathbb{R}^{N}\rightarrow \left[ 0,\infty \right) $ are
locally Hölder continuous functions of exponent $\alpha \in \left(
0,1\right) $;

(C1)\quad $f_{j}:\left[ 0,\infty \right) ^{m}\rightarrow \left[ 0,\infty
\right) $ are continuously differentiable in each variable, $f_{j}\left(
0,...,0\right) =0$, $f_{j}\left( s_{1},...,s_{m}\right) >0$ for $s_{i}>0$ ($%
i=1,...,m$);

(C2)\quad $f_{j}$ are increasing on $\left[ 0,\infty \right) ^{m}$ in each
variable;

(C3)\quad $\int_{1}^{\infty }[F\left( s\right) ]^{-1/p}ds=\infty $ \ \ ($%
F\left( s\right) =\int_{0}^{s}\underset{i=1}{\overset{m}{\Sigma }}%
f_{i}\left( t,...,t\right) dt$).

We are concerned with the following system of quasilinear elliptic partial
differential equations in all of space 
\begin{equation}
\left\{ 
\begin{array}{l}
\Delta _{p}u_{1}\left( x\right) =a_{1}\left( x\right) f_{1}\left(
u_{1}\left( x\right) ,...,u_{m}\left( x\right) \right) \text{ }, \\ 
... \\ 
\Delta _{p}u_{m}\left( x\right) =a_{m}\left( x\right) f_{m}\left(
u_{1}\left( x\right) ,...,u_{m}\left( x\right) \right) \text{ },%
\end{array}%
\right. \text{ }x\in \mathbb{R}^{N}  \label{11}
\end{equation}%
where $N-1\geq p>1$, $\Delta _{p}$ stands for the p-Laplacian operator
defined by $\Delta _{p}u:=\func{div}\left( \left\vert \nabla u\right\vert
^{p-2}\nabla u\right) ,1<p<\infty $. 

Problems related to (\ref{11}) have received an increased interest in the
past several decades (see \cite{BZ,CD,NK,K,LAIR,ME,MA,DY,S,Y} with their
references). On the other hand, it is well known that the form (\ref{11}) is
very simple and easy to be written, many important nonlinear partial
differential equations arising from several areas of mathematics and various
sciences including physics, the generalized reaction diffusion theory and
chemistry sciences take this form. Particularly, one of the most important
classes of such partial differential equations are the \textit{%
time-independent Schrödinger equation in quantum mechanics }%
\begin{equation}
(h^{2}/2p)\Delta u=(V-E)u  \label{s}
\end{equation}%
where $h$ is the Plank constant, $p$ is the mass of a particle moving under
the action of a force field described by the potential $V$ whose wave
function is $u$ and the quantity $E$ is the total energy of the particle,
problems which falls into the class of equations discussed here. \ The
equation (\ref{s}) was invented at a time when electrons, protons and
neutrons were considered to be the elementary particles (see \cite{GG} for
details and some new device ).

The modern structure of the nonlinear Schrödinger equation is much more
complicated. Max Born says: "\textit{Who among us has not written the words
`Schrödinger equation' or `Schrödinger function' countless times? The next
generation will probably do the same, and keep his name alive}" which is
true and in our case even we will refer to more general problem (\ref{11}).

Our objective in the present work, in short, is to complete the principal
results of \cite{LI} and \cite{CD2} and other associated works (see for
example, \cite{BZ,NK,KN,LAIR,LAIR2,LI,MA,DY,Y} and references therein)
showing non-existence and existence of solutions for the similar problem (%
\ref{11}).

Among the obtained results, two of them seem to be worth stressing. The
first one is the problem of existence of a non-radially symmetric bounded
solution to (\ref{11}). The second one is to give a necessary and a
sufficient condition for a positive radial solution of (\ref{11}) to be
large.

We summarize in the next theorems the main objectives of the paper:

\begin{theorem}
\label{1}Assume that\textit{\ }$a_{j}$\textit{\ (}$j=1,...,m$\textit{)
satisfy (A1), }$f_{j}$\textit{\ satisfy (C1)-(C3) and that there exists a
positive number }$\varepsilon $\textit{\ such that }%
\begin{equation}
\int_{0}^{\infty }t^{1+\varepsilon }\left( \underset{j=1}{\overset{m}{\Sigma 
}}\varphi _{j}\left( t\right) \right) ^{2/p}dt<\infty \text{ where }\varphi
_{j}\left( t\right) =\max_{\left\vert x\right\vert =t}a_{j}\left( x\right)
\label{5}
\end{equation}%
\textit{and } $r^{\frac{p\left( N-1\right) }{p-1}}\underset{j=1}{\overset{m}{%
\Sigma }}\varphi _{j}\left( r\right) $\textit{\ is nondecreasing for large }$%
r$\textit{\ then system (\ref{11}) has a nonnegative nontrivial bounded
solution on }$\mathbb{R}^{N}$\textit{. If, on the other hand, }$a_{j}$%
\textit{\ satisfy}%
\begin{equation}
\int_{0}^{\infty }t^{1/\left( p-1\right) }\underset{j=1}{\overset{m}{\Sigma }%
}\psi _{j}^{1/\left( p-1\right) }\left( t\right) dt=\infty \text{ where }%
\psi _{j}\left( t\right) =\min_{\left\vert x\right\vert =t}a_{j}\left(
x\right) ,  \label{5b}
\end{equation}%
\textit{then system (\ref{11}) has no nonnegative nontrivial entire bounded
radial solution on }$\mathbb{R}^{N}$\textit{.}
\end{theorem}

\begin{theorem}
\label{2} \textit{Suppose that }$a_{j}$\textit{\ }$:\left[ 0,\infty \right)
\rightarrow \left[ 0,\infty \right) $\textit{\ (}$j=1,...,m$) \textit{are
continuous spherically symmetric functions (i.e. }$a_{j}\left( x\right)
=a_{j}\left( \left\vert x\right\vert \right) $\textit{). If }$f_{j}$\textit{%
\ satisfy (C1)-(C3) then the problem (\ref{11}) has a non-negative
non-trivial entire radial solution. Suppose furthermore that } $r^{\frac{%
p\left( N-1\right) }{p-1}}\underset{j=1}{\overset{m}{\Sigma }}a_{j}\left(
r\right) $\textit{\ are nondecreasing for large }$r$\textit{. If, in
addition }$a_{j}$\textit{\ satisfies }%
\begin{equation}
\int_{0}^{\infty }\left( \frac{1}{t^{N-1}}\int_{0}^{t}s^{N-1}a_{j}\left(
s\right) ds\right) ^{1/\left( p-1\right) }dt=\infty \text{ for all }j=1,...,m%
\text{ }  \label{12}
\end{equation}%
\textit{then any nonnegative nontrivial solution }$\left(
u_{1},...,u_{m}\right) $\textit{\ of (\ref{11}) is large. Conversely, if (%
\ref{11}) has a nonnegative entire large solution, then }$a_{j}$\textit{\
satisfy}%
\begin{equation}
\int_{0}^{\infty }r^{1+\varepsilon }\left( \underset{j=1}{\overset{d}{\Sigma 
}}a_{j}\left( r\right) \right) ^{2/p}dr=\infty ,  \label{13}
\end{equation}%
\textit{for every }$\varepsilon >0$\textit{.}
\end{theorem}

We remark that in the case $p\neq 2$ the above results are new even for the
situation $m=1$ analyzed in the work \cite{Y}, where the authors have a
problem in the proof when $N-1=p$ which is solved here.

The contents of the paper are organized as follows: in Section \ref{s1} we
establish a preliminary result. Section \ref{ss1} deals with the proof of
Theorem \ref{1}, Section \ref{sss1} is devoted to proving Theorems \ref{2}.

\section{Preliminary results \label{s1}}

For the statement of the next result we need some additional definitions.

\begin{definition}
A function $\left( w_{1},...,w_{m}\right) \in \left[ C_{loc}^{1,\alpha
}\left( \mathbb{R}^{N}\right) \right] ^{m}$ ($\alpha \in \left( 0,1\right) $%
) is called a lower solution of the problem (\ref{11}) if%
\begin{eqnarray*}
\Delta _{p}w_{1}\left( x\right) &\geq &a_{1}\left( x\right) f_{1}\left(
w_{1}\left( x\right) ,...,w_{m}\left( x\right) \right) \text{ for }x\in 
\mathbb{R}^{N}\text{ }, \\
&&... \\
\Delta _{p}w_{m}\left( x\right) &\geq &a_{m}\left( x\right) f_{m}\left(
w_{1}\left( x\right) ,...,w_{m}\left( x\right) \right) \text{ for }x\in 
\mathbb{R}^{N}.
\end{eqnarray*}%
$\ $
\end{definition}

\begin{definition}
We will say that $\left( v_{1},...,v_{m}\right) \in \left[ C_{loc}^{1,\alpha
}\left( \mathbb{R}^{N}\right) \right] ^{m}$ ($\alpha \in \left( 0,1\right) $%
) is an upper solution of the problem (\ref{11}) if%
\begin{eqnarray*}
\Delta _{p}v_{1}\left( x\right) &\leq &a_{1}\left( x\right) f_{1}\left(
v_{1}\left( x\right) ,...,v_{m}\left( x\right) \right) \text{ for }x\in 
\mathbb{R}^{N}\text{ }, \\
&&... \\
\Delta _{p}v_{m}\left( x\right) &\leq &a_{m}\left( x\right) f_{m}\left(
v_{1}\left( x\right) ,...,v_{m}\left( x\right) \right) \text{ for }x\in 
\mathbb{R}^{N}\text{ }.
\end{eqnarray*}
\end{definition}

The proof of Theorem 1.1 is based on the results below, which can be proved
as in \cite[Theorem 5.1, pp. 146]{NK}, \cite[Lemma 1, pp. 15]{QY} or \cite%
{PAO} in the same time with \cite{ME}.

\begin{lemma}
\label{lu}Make the same assumptions on $a_{j}$ and \ $f_{j}$ ($j=1,...,m$)
as in Theorem \ref{1}. If the problem (\ref{11}) has a pair of upper and
lower bounded solutions $\left( v_{1},...,v_{m}\right) $ and $\left(
w_{1},...,w_{m}\right) $ fulfilling $w_{j}(x)\leq v_{j}(x),$ ($j=1,...,m$)$%
,\forall x\in \mathbb{R}^{N}$ then there exists a bounded function $\left(
u_{1},...,u_{m}\right) $ belonging to $\left[ C_{loc}^{1,\alpha }\left( 
\mathbb{R}^{N}\right) \right] ^{m}$ ($\alpha \in \left( 0,1\right) $) with%
\begin{equation*}
w_{j}(x)\leq u_{j}(x)\leq v_{j}(x)\text{, }\forall x\in \mathbb{R}^{N}
\end{equation*}%
and satisfying (\ref{11}).
\end{lemma}

\section{Proof of the Theorem \protect\ref{1} \label{ss1}}

The proof is inspired by the corresponding ones for the $p=2$ cases in \cite%
{CD2} with some new ideas. Assume that (\ref{5}) holds. We use the method of
lower and upper solutions for the problem (\ref{11}). We look for an upper
solution $\left( v_{1},...,v_{m}\right) $ and a lower solution $\left(
w_{1},...,w_{m}\right) $. To find a positive lower solution, we observe that
an arbitrary positive solution $w_{i}$ ($i=1,...,m$)\ to the following
auxiliary system%
\begin{equation}
\left\{ 
\begin{array}{c}
\Delta _{p}w_{1}\left( r\right) =\varphi _{1}\left( r\right) f_{i}\left(
w_{1},...,w_{m}\right) \text{ for }r:=\left\vert x\right\vert \text{, }x\in 
\mathbb{R}^{N}\text{ }, \\ 
... \\ 
\Delta _{p}w_{m}\left( r\right) =\varphi _{m}\left( r\right) f_{m}\left(
w_{1},...,w_{m}\right) \text{ for }r:=\left\vert x\right\vert \text{, }x\in 
\mathbb{R}^{N}\text{ },%
\end{array}%
\right. \text{ }  \label{6}
\end{equation}%
is the best candidate. We shall only study the radial solutions of (\ref{6}%
), hence always write (\ref{6}) in the following radial version: 
\begin{eqnarray}
\left( p-1\right) w_{1}^{\prime }\left( r\right) ^{p-2}w_{1}^{\prime \prime
}+\frac{N-1}{r}w_{1}^{\prime }\left( r\right) ^{p-1} &=&\varphi _{1}\left(
r\right) f_{1}\left( w_{1}\left( r\right) ,...,w_{m}\left( r\right) \right) ,
\notag \\
&&....  \label{66} \\
\left( p-1\right) w_{m}^{\prime }\left( r\right) ^{p-2}w_{m}^{\prime \prime
}+\frac{N-1}{r}w_{m}^{\prime }\left( r\right) ^{p-1} &=&\varphi _{m}\left(
r\right) f_{m}\left( w_{1}\left( r\right) ,...,w_{m}\left( r\right) \right) .
\notag
\end{eqnarray}%
First we see that radial solutions of (\ref{66}) are any positive solutions $%
\left( w_{1},...,w_{m}\right) $ of the integral equations%
\begin{eqnarray*}
w_{1}\left( r\right)  &=&\frac{1}{m}+\int_{0}^{r}\left( \frac{1}{t^{N-1}}%
\int_{0}^{t}s^{N-1}\varphi _{1}\left( s\right) f_{1}\left( w_{1}\left(
s\right) ,...,w_{m}\left( s\right) \right) ds\right) ^{1/\left( p-1\right)
}dt, \\
&&... \\
w_{m}\left( r\right)  &=&\frac{1}{m}+\int_{0}^{r}\left( \frac{1}{t^{N-1}}%
\int_{0}^{t}s^{N-1}\varphi _{m}\left( s\right) f_{m}\left( w_{1}\left(
s\right) ,...,w_{m}\left( s\right) \right) ds\right) ^{1/\left( p-1\right)
}dt.
\end{eqnarray*}%
To establish a solution to this system, we use successive approximation.
Define, recursively, sequences $\left\{ w_{i}^{k}\right\} _{i=\overline{1,m}%
}^{k\geq 1}$ on $\left[ 0,\infty \right) $ by%
\begin{equation*}
\left\{ 
\begin{array}{l}
w_{1}^{0}=...=w_{m}^{0}=\frac{1}{m},\text{ }r\geq 0, \\ 
w_{1}^{k}\left( r\right) =\frac{1}{m}+\int_{0}^{r}\left( \frac{1}{t^{N-1}}%
\int_{0}^{t}s^{N-1}\varphi _{1}\left( s\right) f_{1}\left( w_{1}^{k-1}\left(
s\right) ,...,w_{m}^{k-1}\left( s\right) \right) ds\right) ^{1/\left(
p-1\right) }dt \\ 
... \\ 
w_{m}^{k}\left( r\right) =\frac{1}{m}+\int_{0}^{r}\left( \frac{1}{t^{N-1}}%
\int_{0}^{t}s^{N-1}\varphi _{m}\left( s\right) f_{m}\left( w_{1}^{k-1}\left(
s\right) ,...,w_{m}^{k-1}\left( s\right) \right) ds\right) ^{1/\left(
p-1\right) }dt.%
\end{array}%
\right. 
\end{equation*}%
We remark that, for all $r\geq 0,$ $i=1,...,m$ and $k\in N$ 
\begin{equation*}
w_{i}^{k}\left( r\right) \geq \frac{1}{m}\text{,}
\end{equation*}%
and that $\left\{ w_{i}^{k}\right\} _{i=1,...,m}^{k\geq 1}$ are
non-decreasing sequence on $\left[ 0,\infty \right) $.

We note that $\left\{ w_{i}^{k}\right\} _{i=1,...,m}^{k\geq 1}$ satisfy%
\begin{eqnarray}
\left( p-1\right) \left[ \left( w_{1}^{k}\right) ^{\prime }\right]
^{p-2}\left( w_{1}^{k}\right) ^{\prime \prime }+\frac{N-1}{r}\left[ \left(
w_{1}^{k}\right) ^{\prime }\right] ^{p-1} &=&\varphi _{1}\left( r\right)
f_{1}\left( w_{1}^{k-1}\left( r\right) ,...,w_{m}^{k-1}\left( r\right)
\right) ,  \notag \\
&&...  \label{sis1} \\
\left( p-1\right) \left[ \left( w_{m}^{k}\right) ^{\prime }\right]
^{p-2}\left( w_{m}^{k}\right) ^{\prime \prime }+\frac{N-1}{r}\left[ \left(
w_{m}^{k}\right) ^{\prime }\right] ^{p-1} &=&\varphi _{m}\left( r\right)
f_{m}\left( w_{1}^{k-1}\left( r\right) ,...,w_{m}^{k-1}\left( r\right)
\right) .  \notag
\end{eqnarray}%
Using the monotonicity of $\left\{ w_{i}^{k}\right\} _{i=1,...,m}^{k\geq 1}$
yields%
\begin{eqnarray}
\varphi _{1}\left( r\right) f_{1}\left( w_{1}^{k-1}\left( r\right)
,...,w_{m}^{k-1}\left( r\right) \right) &\leq &\varphi _{1}\left( r\right)
f_{1}\left( w_{1}^{k},...,w_{m}^{k}\right) \leq \varphi _{1}\left( r\right) 
\overset{m}{\underset{i=1}{\Sigma }}f_{i}\left( \overset{m}{\underset{i=1}{%
\Sigma }}w_{i}^{k},...,\overset{m}{\underset{i=1}{\Sigma }}w_{i}^{k}\right) ,
\notag \\
&&...  \label{8} \\
\varphi _{m}\left( r\right) f_{m}\left( w_{1}^{k-1}\left( r\right)
,...,w_{m}^{k-1}\left( r\right) \right) &\leq &\varphi _{m}\left( r\right)
f_{i}\left( w_{1}^{k},...,w_{m}^{k}\right) \leq \varphi _{m}\left( r\right) 
\overset{m}{\underset{i=1}{\Sigma }}f_{i}\left( \overset{m}{\underset{i=1}{%
\Sigma }}w_{i}^{k},...,\overset{m}{\underset{i=1}{\Sigma }}w_{i}^{k}\right) ,
\notag
\end{eqnarray}%
and, so%
\begin{eqnarray}
\left( p-1\right) \left[ \left( w_{1}^{k}\left( r\right) \right) ^{\prime }%
\right] ^{p-1}\left( w_{1}^{k}\right) ^{\prime \prime }+\frac{N-1}{r}\left[
\left( w_{1}^{k}\left( r\right) \right) ^{\prime }\right] ^{p} &\leq
&\varphi _{1}\left( r\right) \overset{m}{\underset{i=1}{\Sigma }}f_{i}\left( 
\overset{m}{\underset{i=1}{\Sigma }}w_{i}^{k},..,\overset{m}{\underset{i=1}{%
\Sigma }}w_{i}^{k}\right) \left( w_{1}^{k}\left( r\right) \right) ^{\prime },
\notag \\
&&...  \label{88} \\
\left( p-1\right) \left[ \left( w_{m}^{k}\left( r\right) \right) ^{\prime }%
\right] ^{p-1}\left( w_{m}^{k}\right) ^{\prime \prime }+\frac{N-1}{r}\left[
\left( w_{m}^{k}\left( r\right) \right) ^{\prime }\right] ^{p} &\leq
&\varphi _{m}\left( r\right) \overset{m}{\underset{i=1}{\Sigma }}f_{i}\left( 
\overset{m}{\underset{i=1}{\Sigma }}w_{i}^{k},..,\overset{m}{\underset{i=1}{%
\Sigma }}w_{i}^{k}\right) \left( w_{m}^{k}\left( r\right) \right) ^{\prime },
\notag
\end{eqnarray}%
which implies that%
\begin{eqnarray}
\left( p-1\right) \left[ \left( w_{1}^{k}\left( r\right) \right) ^{\prime }%
\right] ^{p-1}\left( w_{1}^{k}\right) ^{\prime \prime }+\frac{N-1}{r}\left[
\left( w_{1}^{k}\left( r\right) \right) ^{\prime }\right] ^{p} &\leq
&\varphi _{1}\left( r\right) \overset{m}{\underset{i=1}{\Sigma }}f_{i}\left( 
\overset{m}{\underset{i=1}{\Sigma }}w_{i}^{k},..,\overset{m}{\underset{i=1}{%
\Sigma }}w_{i}^{k}\right) \left( \overset{m}{\underset{i=1}{\Sigma }}%
w_{i}^{k}\left( r\right) \right) ^{\prime },  \notag \\
&&...  \label{sum} \\
\left( p-1\right) \left[ \left( w_{m}^{k}\left( r\right) \right) ^{\prime }%
\right] ^{p-1}\left( w_{m}^{k}\right) ^{\prime \prime }+\frac{N-1}{r}\left[
\left( w_{m}^{k}\left( r\right) \right) ^{\prime }\right] ^{p} &\leq
&\varphi _{m}\left( r\right) \overset{m}{\underset{i=1}{\Sigma }}f_{i}\left( 
\overset{m}{\underset{i=1}{\Sigma }}w_{i}^{k},..,\overset{m}{\underset{i=1}{%
\Sigma }}w_{i}^{k}\right) \left( \overset{m}{\underset{i=1}{\Sigma }}%
w_{i}^{k}\left( r\right) \right) ^{\prime }.  \notag
\end{eqnarray}%
We choose $R>0$ so that $r^{\frac{p\left( N-1\right) }{p-1}}\underset{j=1}{%
\overset{m}{\Sigma }}\varphi _{j}\left( r\right) $ are non-decreasing for $%
r\geq R$. First we prove that $w_{i}^{k}\left( R\right) $ and $\left(
w_{i}^{k}\left( R\right) \right) ^{\prime }$, both of which are nonnegative,
are bounded above independent of $k$. To do this, let 
\begin{equation*}
\phi _{i}^{R}=\max \{\varphi _{i}\left( r\right) :0\leq r\leq R\}\text{, }%
i=1,...,m\text{.}
\end{equation*}%
Using this and the fact that $\left( w_{i}^{k}\right) ^{\prime }\geq 0$ ($%
i=1,...,m$), we note that (\ref{sum}) and (\ref{88}) yields%
\begin{eqnarray*}
\left( p-1\right) \left[ \left( w_{1}^{k}\left( r\right) \right) ^{\prime }%
\right] ^{p-1}\left( w_{1}^{k}\right) ^{\prime \prime }+\frac{N-1}{r}\left[
\left( w_{1}^{k}\left( r\right) \right) ^{\prime }\right] ^{p} &\leq &\phi
_{1}^{R}\overset{m}{\underset{i=1}{\Sigma }}f_{i}\left( \overset{m}{\underset%
{i=1}{\Sigma }}w_{i}^{k},...,\overset{m}{\underset{i=1}{\Sigma }}%
w_{i}^{k}\right) \left( \overset{m}{\underset{i=1}{\Sigma }}w_{i}^{k}\left(
r\right) \right) ^{\prime } \\
&&... \\
\left( p-1\right) \left[ \left( w_{m}^{k}\left( r\right) \right) ^{\prime }%
\right] ^{p-1}\left( w_{m}^{k}\right) ^{\prime \prime }+\frac{N-1}{r}\left[
\left( w_{m}^{k}\left( r\right) \right) ^{\prime }\right] ^{p} &\leq &\phi
_{m}^{R}\overset{m}{\underset{i=1}{\Sigma }}f_{i}\left( \overset{m}{\underset%
{i=1}{\Sigma }}w_{i}^{k},...,\overset{m}{\underset{i=1}{\Sigma }}%
w_{i}^{k}\right) \left( \overset{m}{\underset{i=1}{\Sigma }}w_{i}^{k}\left(
r\right) \right) ^{\prime }
\end{eqnarray*}%
which implies%
\begin{equation*}
\left\{ \underset{i=1}{\overset{m}{\Sigma }}\left[ \left( w_{i}^{k}\left(
r\right) \right) ^{\prime }\right] ^{p}\right\} ^{\prime }\leq \frac{p}{p-1}%
\underset{i=1}{\overset{m}{\Sigma }}\phi _{i}^{R}\overset{m}{\underset{i=1}{%
\Sigma }}f_{i}\left( \overset{m}{\underset{i=1}{\Sigma }}w_{i}^{k}\left(
r\right) ,...,\overset{m}{\underset{i=1}{\Sigma }}w_{i}^{k}\left( r\right)
\right) \left( \underset{i=1}{\overset{m}{\Sigma }}w_{i}^{k}\left( r\right)
\right) ^{\prime }.
\end{equation*}%
Integrate this equation from $0$ to $r$. We obtain%
\begin{equation}
\underset{i=1}{\overset{m}{\Sigma }}\left[ \left( w_{i}^{k}\left( r\right)
\right) ^{\prime }\right] ^{p}\leq \frac{p}{p-1}\underset{i=1}{\overset{m}{%
\Sigma }}\phi _{i}^{R}\int_{1}^{\overset{m}{\underset{i=1}{\Sigma }}%
w_{i}^{k}\left( r\right) }\overset{m}{\underset{i=1}{\Sigma }}f_{i}\left(
s,...,s\right) ds\text{, }0\leq r\leq R.  \label{ineq2}
\end{equation}%
Since $p>1$ we know that%
\begin{equation}
\left( a_{1}+...+a_{m}\right) ^{p}\leq m^{p-1}\left(
a_{1}^{p}+...+a_{m}^{p}\right)  \label{ineq}
\end{equation}%
for any non-negative constants $a_{i}$ ($i=1,...,m$). Using this inequality
in (\ref{ineq2}) we have%
\begin{equation*}
m^{1-p}\left[ \underset{i=1}{\overset{m}{\Sigma }}\left( w_{i}^{k}\left(
r\right) \right) ^{\prime }\right] ^{p}\leq \frac{p}{p-1}\underset{i=1}{%
\overset{m}{\Sigma }}\phi _{i}^{R}\int_{1}^{\overset{m}{\underset{i=1}{%
\Sigma }}w_{i}^{k}\left( r\right) }\overset{m}{\underset{i=1}{\Sigma }}%
f_{i}\left( s,...,s\right) ds\text{, }0\leq r\leq R,
\end{equation*}%
which yields%
\begin{equation}
\left( \underset{i=1}{\overset{m}{\Sigma }}w_{i}^{k}\left( r\right) \right)
^{\prime }\leq \sqrt[p]{\frac{pm^{p-1}}{p-1}\underset{i=1}{\overset{m}{%
\Sigma }}\phi _{i}^{R}}\left( \int_{1}^{\overset{m}{\underset{i=1}{\Sigma }}%
w_{i}^{k}\left( r\right) }\overset{m}{\underset{i=1}{\Sigma }}f_{i}\left(
s,...,s\right) ds\right) ^{1/p}\text{, }0\leq r\leq R.  \label{9}
\end{equation}%
Integrating the above equation between $0$ and $R$, we have 
\begin{equation*}
\int_{1}^{\overset{m}{\underset{i=1}{\Sigma }}w_{i}^{k}\left( R\right) }%
\left[ \int_{1}^{t}\overset{m}{\underset{i=1}{\Sigma }}f_{i}\left(
s,...,s\right) ds\right] ^{-1/p}dt\leq \sqrt[p]{\frac{pm^{p-1}}{p-1}\left( 
\underset{i=1}{\overset{m}{\Sigma }}\phi _{i}^{R}\right) }R\text{.}
\end{equation*}%
By the assumption C3), we now conclude that $\overset{m}{\underset{i=1}{%
\Sigma }}w_{i}^{k}\left( R\right) $ is bounded above independent of $k$ and
using this fact in (\ref{9}) shows that the same is true of $\left( \overset{%
m}{\underset{i=1}{\Sigma }}w_{i}^{k}\left( R\right) \right) ^{\prime }$.
Thus, the sequences $w_{i}^{k}\left( R\right) $ and $\left( w_{i}^{k}\left(
R\right) \right) ^{\prime }$ are bounded above independent of $k$. Now let
us verify that the non-decreasing sequences $w_{i}^{k}$ is bounded for all $%
r\geq 0$ and all $k$. Multiplying (\ref{sum}) by $\frac{p}{p-1}r^{\frac{%
p\left( N-1\right) }{p-1}}$ and summing we have%
\begin{equation*}
\left\{ r^{\frac{p\left( N-1\right) }{p-1}}\overset{m}{\underset{i=1}{\Sigma 
}}\left[ \left( w_{i}^{k}\right) ^{\prime }\right] ^{p}\right\} ^{\prime
}\leq \text{ }\frac{pr^{\frac{p\left( N-1\right) }{p-1}}}{p-1}\overset{m}{%
\underset{i=1}{\Sigma }}\varphi _{i}\left( r\right) \overset{m}{\underset{i=1%
}{\Sigma }}f_{i}\left( \overset{m}{\underset{i=1}{\Sigma }}w_{i}^{k},...,%
\overset{m}{\underset{i=1}{\Sigma }}w_{i}^{k}\right) \left( \overset{m}{%
\underset{i=1}{\Sigma }}w_{i}^{k}\right) ^{\prime }.
\end{equation*}%
and integrating gives%
\begin{eqnarray}
&&\int_{R}^{r}\left\{ s^{\frac{p\left( N-1\right) }{p-1}}\overset{m}{%
\underset{i=1}{\Sigma }}\left[ \left( w_{i}^{k}\left( s\right) \right)
^{\prime }\right] ^{p}\right\} ^{\prime }ds  \notag \\
&\leq &\int_{R}^{r}\frac{pm^{p-1}}{p-1}s^{\frac{p\left( N-1\right) }{p-1}}%
\overset{m}{\underset{i=1}{\Sigma }}\varphi _{i}\left( s\right) \overset{m}{%
\underset{i=1}{\Sigma }}f_{i}\left( \overset{m}{\underset{i=1}{\Sigma }}%
w_{i}^{k}\left( s\right) ,...,\overset{m}{\underset{i=1}{\Sigma }}%
w_{i}^{k}\left( s\right) \right) \left( \overset{m}{\underset{i=1}{\Sigma }}%
w_{i}^{k}\right) ^{\prime }ds.  \label{s2}
\end{eqnarray}%
Hence (\ref{ineq}) in (\ref{s2}) gives 
\begin{eqnarray*}
&&r^{\frac{p\left( N-1\right) }{p-1}}m^{1-p}\left[ \left( \overset{m}{%
\underset{i=1}{\Sigma }}w_{i}^{k}\left( r\right) \right) ^{\prime }\right]
^{p}-R^{\frac{p\left( N-1\right) }{p-1}}\overset{m}{\underset{i=1}{\Sigma }}%
\left[ \left( w_{i}^{k}\left( R\right) \right) ^{\prime }\right] ^{p} \\
&\leq &\int_{R}^{r}\frac{p}{p-1}s^{\frac{p\left( N-1\right) }{p-1}}\overset{m%
}{\underset{i=1}{\Sigma }}\varphi _{i}\left( s\right) \overset{m}{\underset{%
i=1}{\Sigma }}f_{i}\left( \overset{m}{\underset{i=1}{\Sigma }}%
w_{i}^{k}\left( s\right) ,...,\overset{m}{\underset{i=1}{\Sigma }}%
w_{i}^{k}\left( s\right) \right) \left( \overset{m}{\underset{i=1}{\Sigma }}%
w_{i}^{k}\right) ^{\prime }ds
\end{eqnarray*}%
and thus%
\begin{eqnarray*}
&&r^{\frac{p\left( N-1\right) }{p-1}}\left[ \left( \overset{m}{\underset{i=1}%
{\Sigma }}w_{i}^{k}\left( r\right) \right) ^{\prime }\right] ^{p}\leq R^{%
\frac{p\left( N-1\right) }{p-1}}m^{p-1}\overset{m}{\underset{i=1}{\Sigma }}%
\left[ \left( w_{i}^{k}\left( R\right) \right) ^{\prime }\right] ^{p}+ \\
&&+\int_{R}^{r}\frac{pm^{p-1}s^{\frac{p\left( N-1\right) }{p-1}}}{p-1}%
\overset{m}{\underset{i=1}{\Sigma }}\varphi _{i}\left( s\right) \overset{m}{%
\underset{i=1}{\Sigma }}f_{i}\left( \overset{m}{\underset{i=1}{\Sigma }}%
w_{i}^{k}\left( s\right) ,...,\overset{m}{\underset{i=1}{\Sigma }}%
w_{i}^{k}\left( s\right) \right) \left( \overset{m}{\underset{i=1}{\Sigma }}%
w_{i}^{k}\right) ^{\prime }ds.
\end{eqnarray*}%
for $r\geq R$. Noting that, by the monotonicity of $s^{\frac{p\left(
N-1\right) }{p-1}}\overset{m}{\underset{i=1}{\Sigma }}\varphi _{i}\left(
s\right) $ for $r\geq s\geq R$, we get 
\begin{equation*}
r^{^{\frac{p\left( N-1\right) }{p-1}}}\left[ \left( \overset{m}{\underset{i=1%
}{\Sigma }}w_{i}^{k}\left( r\right) \right) ^{\prime }\right] ^{p}\leq C+%
\frac{pm^{p-1}}{p-1}r^{^{\frac{p\left( N-1\right) }{p-1}}}\overset{m}{%
\underset{i=1}{\Sigma }}\varphi _{i}\left( r\right) F\left( \overset{m}{%
\underset{i=1}{\Sigma }}w_{i}^{k}\left( r\right) \right) ,
\end{equation*}%
where 
\begin{equation*}
C=R^{^{\frac{p\left( N-1\right) }{p-1}}}m^{p-1}\overset{m}{\underset{i=1}{%
\Sigma }}\left[ \left( w_{i}^{k}\left( R\right) \right) ^{\prime }\right]
^{p},
\end{equation*}%
which yields%
\begin{equation}
\left( \overset{m}{\underset{i=1}{\Sigma }}w_{i}^{k}\right) ^{\prime }\leq %
\left[ Cr^{^{\frac{p\left( 1-N\right) }{p-1}}}+\frac{pm^{p-1}}{p-1}\overset{m%
}{\underset{i=1}{\Sigma }}\varphi _{i}\left( r\right) F\left( \overset{m}{%
\underset{i=1}{\Sigma }}w_{i}^{k}\left( r\right) \right) \right] ^{1/p}.
\label{100}
\end{equation}%
Since $\left( 1/p\right) <1$ we know that%
\begin{equation*}
\left( b_{1}+b_{2}\right) ^{1/p}\leq b_{1}^{1/p}+b_{2}^{1/p}
\end{equation*}%
for any non-negative constants $b_{i}$ ($i=1,...,m$). Therefore, by applying
this inequality in (\ref{100}) we get%
\begin{equation*}
\left( \overset{m}{\underset{i=1}{\Sigma }}w_{i}^{k}\right) ^{\prime }\leq 
\sqrt[p]{C}r^{\left( 1-N\right) /(p-1)}+\sqrt[p]{\frac{pm^{p-1}}{p-1}\overset%
{m}{\underset{i=1}{\Sigma }}\varphi _{i}\left( r\right) }\left[ F\left( 
\overset{m}{\underset{i=1}{\Sigma }}w_{i}^{k}\left( r\right) \right) \right]
^{1/p}.
\end{equation*}%
Integrating the above inequality, we get%
\begin{equation}
\frac{m}{mr}\int_{\overset{m}{\underset{i=1}{\Sigma }}w_{i}^{k}\left(
R\right) }^{\overset{m}{\underset{i=1}{\Sigma }}w_{i}^{k}\left( r\right) }%
\left[ F\left( t\right) \right] ^{-1/p}dt\leq \sqrt[p]{C}r^{\left(
1-N\right) /\left( p-1\right) }\left[ F\left( \overset{m}{\underset{i=1}{%
\Sigma }}w_{i}^{k}\left( r\right) \right) \right] ^{-1/p}+\left( \frac{%
pm^{p-1}}{p-1}\overset{m}{\underset{i=1}{\Sigma }}\varphi _{i}\left(
r\right) \right) ^{1/p}.  \label{sis3}
\end{equation}%
Integrating (\ref{sis3}) and using the fact that%
\begin{equation*}
\left( \overset{m}{\underset{i=1}{\Sigma }}\varphi _{i}\left( s\right)
\right) ^{1/p}=\left( s^{p\left( 1+\varepsilon \right) /2}\overset{m}{%
\underset{i=1}{\Sigma }}\varphi _{i}\left( s\right) s^{-p\left(
1+\varepsilon \right) /2}\right) ^{1/p}\leq \left( \frac{1}{2}\right) ^{1/p}%
\left[ s^{1+\varepsilon }\left( \overset{m}{\underset{i=1}{\Sigma }}\varphi
_{i}\left( r\right) \right) ^{2/p}+s^{-1-\varepsilon }\right]
\end{equation*}%
for each $\varepsilon >0$, we have%
\begin{eqnarray}
\int_{\overset{m}{\underset{i=1}{\Sigma }}w_{i}^{k}\left( R\right) }^{%
\overset{m}{\underset{i=1}{\Sigma }}w_{i}^{k}\left( r\right) }\left[ F\left(
t\right) \right] ^{-1/p}dt &\leq &\sqrt[p]{C}\int_{R}^{r}t^{\frac{1-N}{p-1}}%
\left[ F\left( \overset{m}{\underset{i=1}{\Sigma }}w_{i}^{k}\left( t\right)
\right) \right] ^{-1/p}dt  \notag \\
&&+\left( \frac{1}{2}\right) ^{1/p}\sqrt[p]{\frac{pm^{p-1}}{p-1}}\left[
\int_{R}^{r}t^{1+\varepsilon }\left( \overset{m}{\underset{i=1}{\Sigma }}%
\varphi _{i}\left( t\right) \right) ^{2/p}dt+\int_{R}^{r}t^{-1-\varepsilon
}dt\right]  \notag \\
&\leq &\sqrt[p]{C}\left[ F\left( \overset{m}{\underset{i=1}{\Sigma }}%
w_{i}^{k}\left( R\right) \right) \right] ^{-1/p}\frac{p-1}{p-N}R^{\frac{p-N}{%
p-1}}  \notag \\
&&+\left( \frac{1}{2}\right) ^{1/p}\sqrt[p]{\frac{pm^{p-1}}{p-1}}\left[
\int_{R}^{r}t^{1+\varepsilon }\left( \overset{m}{\underset{i=1}{\Sigma }}%
\varphi _{i}\left( t\right) \right) ^{2/p}dt+\frac{1}{\varepsilon
R^{\varepsilon }}\right] \text{.}  \label{111}
\end{eqnarray}%
Since the right side of this inequality is bounded independent of $k$ (note
that $w_{i}^{k}\left( t\right) \geq 1/m$), so is the left side and hence, in
light of C3), the sequence $\left\{ \overset{m}{\underset{i=1}{\Sigma }}%
w_{i}^{k}\right\} ^{k\geq 1}$ is a bounded sequence and so $\left\{
w_{j}^{k}\right\} _{j=1,...,m}^{k\geq 1}$ are bounded sequence. Thus, for
every $x\in \mathbb{R}^{N}$, it makes sense to define $w_{j}\left(
\left\vert x\right\vert \right) :=\underset{k\rightarrow \infty }{\lim }%
w_{j}^{k}\left( \left\vert x\right\vert \right) $ for all $j=1,...,m$ and so 
$\left( w_{1},...,w_{m}\right) $ is a positive solution of (\ref{6}).

Since, we have found upper bounds for $\left\{ w_{j}\right\} _{j=1,...,m}$
we can let $M$ be the least upper bound of $\overset{m}{\underset{i=1}{%
\Sigma }}w_{i}$ and note that 
\begin{equation*}
M=\lim_{r\rightarrow \infty }\overset{m}{\underset{i=1}{\Sigma }}w_{i}\left(
r\right) \text{.}
\end{equation*}%
Now let $\psi _{i}\left( t\right) =\min_{\left\vert x\right\vert
=t}a_{i}\left( x\right) $ and $v_{i}$ ($i=1,...,m$)\ be the positive
increasing bounded solutions of 
\begin{eqnarray*}
v_{1}\left( r\right)  &=&M+\int_{0}^{r}\left( \frac{1}{t^{N-1}}%
\int_{0}^{t}s^{N-1}\psi _{1}\left( s\right) f_{1}\left( v_{1}\left( s\right)
,...,v_{m}\left( s\right) \right) ds\right) ^{1/p-1}dt,\text{ } \\
&&... \\
v_{m}\left( r\right)  &=&M+\int_{0}^{r}\left( \frac{1}{t^{N-1}}%
\int_{0}^{t}s^{N-1}\psi _{m}\left( s\right) f_{m}\left( v_{1}\left( s\right)
,...,v_{m}\left( s\right) \right) ds\right) ^{1/p-1}dt,
\end{eqnarray*}%
which, of course, satisfies (\ref{6}) with $w_{i}$ replaced with $v_{i}$ and 
$\varphi _{i}$ replaced with $\psi _{i}$. It is also clear that $v_{i}\geq M$
($i=1,...,m$). (The proof of the existence of $v_{i}$ and that it has the
properties mentioned is similar to that given in the proof for $w_{i}$ \ and
is therefore omitted). Thus we have an upper solution $\left(
v_{1},...,v_{m}\right) $ and a lower solution $\left( w_{1},...,w_{m}\right) 
$. Then the standard upper-lower solution principle (see Lemma \ref{lu})
implies that the problem (\ref{11}) has a solution $\left(
u_{1},...,u_{m}\right) $.

We end this section analyzing the non-existence of solutions. For this,
assume that (\ref{5b}) holds. Arguing by contradiction, let us assume that
the system (\ref{11}) has nonnegative nontrivial entire bounded solution $%
\left( u_{1},...,u_{m}\right) $ on $\mathbb{R}^{N}$. \ Assuming $%
M_{i}=\sup_{x\in \mathbb{R}^{N}}u_{i}\left( x\right) $ ($i=1,...,m$) and
knowing that $u_{i}^{\prime }\geq 0$, we get $\lim_{r\rightarrow \infty
}u_{i}\left( r\right) =M_{i}$ ($i=1,...,m$). Thus there exists $R>0$ such
that $u_{i}\geq M_{i}/2$ for $r\geq R$. From conditions of $f_{i}$, it
follows that 
\begin{eqnarray}
\text{ }f_{1}^{1/\left( p-1\right) }\left( u_{1},...,u_{m}\right)  &\geq
&f_{1}^{1/\left( p-1\right) }\left( M_{1}/2,...,M_{m}/2\right) :=c_{0}^{1}%
\text{ }  \notag \\
&&  \label{con} \\
f_{m}^{1/\left( p-1\right) }\left( u_{1},...,u_{m}\right)  &\geq
&f_{m}^{1/\left( p-1\right) }\left( M_{1}/2,...,M_{m}/2\right) :=c_{0}^{m} 
\notag
\end{eqnarray}%
for $r\geq R$. On the other hand%
\begin{eqnarray}
u_{1} &\geq &u_{1}\left( 0\right) +\int_{0}^{r}\left(
t^{1-N}\int_{0}^{t}s^{N-1}\psi _{1}\left( s\right) f_{1}\left( u_{1}\left(
s\right) ,...,u_{m}\left( s\right) \right) ds\right) ^{1/\left( p-1\right)
}dt,  \notag \\
&&...  \label{sis4} \\
u_{m} &\geq &u_{m}\left( 0\right) +\int_{0}^{r}\left(
t^{1-N}\int_{0}^{t}s^{N-1}\psi _{m}\left( s\right) f_{m}\left( u_{1}\left(
s\right) ,...,u_{m}\left( s\right) \right) ds\right) ^{1/\left( p-1\right)
}dt.  \notag
\end{eqnarray}%
Rearranging the terms, and by using the conditions (\ref{con}) in (\ref{sis4}%
) follows%
\begin{equation*}
\overset{m}{\underset{i=1}{\Sigma }}u_{i}\left( r\right) \geq
mc_{1}+c_{2}\left( \frac{1}{N}\right) ^{1/\left( p-1\right)
}\int_{R}^{r}t^{1-N}\overset{m}{\underset{i=1}{\Sigma }}\psi _{i}^{1/\left(
p-1\right) }\left( t\right) dt\rightarrow \infty \text{ as }r\rightarrow
\infty ,\text{ }
\end{equation*}%
where $c_{1}=\min \{u_{1}\left( 0\right) ,...,u_{m}\left( 0\right) \}$ and $%
c_{2}:=\min \{c_{0}^{1},...,c_{0}^{m}\}$ (see also \cite{BZ}). A
contradiction to the boundedness of $\overset{m}{\underset{i=1}{\Sigma }}%
u_{i}\left( r\right) $. This concludes the proof. \rule{5pt}{5pt}

\section{Proof of the Theorem \protect\ref{2} \label{sss1}}

We first notice that from \cite[Theorem 2]{Y} the problem%
\begin{equation}
\Delta _{p}z\left( r\right) =\overset{m}{\underset{i=1}{\Sigma }}a_{i}\left(
r\right) \overset{m}{\underset{i=1}{\Sigma }}f_{i}\left( z\left( r\right)
,...,z\left( r\right) \right) \text{ for }r:=\left\vert x\right\vert \text{, 
}x\in \mathbb{R}^{N}  \label{prob}
\end{equation}%
has a non-negative non-trivial entire solution. Moreover, for each $R>0$,
there exists $c_{R}>0$ such that $z\left( R\right) \leq c_{R}$. Due to the
fact that $z$ is radial, we have%
\begin{equation*}
z\left( r\right) =z\left( 0\right) +\int_{0}^{r}\frac{1}{t^{N-1}}\left(
\int_{0}^{t}s^{N-1}\overset{m}{\underset{i=1}{\Sigma }}a_{i}\left( s\right) 
\overset{m}{\underset{i=1}{\Sigma }}f_{i}\left( z\left( s\right)
,...,z\left( s\right) \right) ds\right) ^{1/\left( p-1\right) }dt\text{ for
all }r\geq 0.
\end{equation*}%
We choose $\beta _{1}\in \left( 0,z\left( 0\right) \right] $. Define the
sequences $\left\{ u_{i}^{k}\right\} _{i=1,...,m}^{k\geq 1}$ on $\left[
0,\infty \right) $ by%
\begin{equation*}
\left\{ 
\begin{array}{l}
u_{1}^{0}=...=u_{m}^{0}=\beta _{1}\text{ for all }r\geq 0 \\ 
u_{1}^{k}\left( r\right) =\beta _{1}+\int_{0}^{r}\left( \frac{1}{t^{N-1}}%
\int_{0}^{t}s^{N-1}a_{1}\left( s\right) f_{1}\left( u_{1}^{k-1}\left(
s\right) ,...,u_{m}^{k-1}\left( s\right) \right) \right) ^{1/\left(
p-1\right) }dsdt, \\ 
... \\ 
u_{m}^{k}\left( r\right) =\beta _{1}+\int_{0}^{r}\left( \frac{1}{t^{N-1}}%
\int_{0}^{t}s^{N-1}a_{m}\left( s\right) f_{m}\left( u_{1}^{k-1}\left(
s\right) ,...,u_{m}^{k-1}\left( s\right) \right) \right) ^{1/\left(
p-1\right) }dsdt.%
\end{array}%
\right. 
\end{equation*}%
With the same arguments as in the proof of Theorem \ref{1} we obtain
that\bigskip\ $\left\{ u_{i}^{k}\right\} _{i=1,...,m}^{k\geq 1}$ are
non-decreasing sequence on $\left[ 0,\infty \right) $. Because $z^{\prime
}\left( r\right) \geq 0$ follows $0<\beta _{1}\leq z\left( 0\right) \leq
z\left( r\right) $ for all $r\geq 0$ and so%
\begin{eqnarray*}
u_{i}^{1}\left( r\right)  &=&\beta _{1}+\int_{0}^{r}\left( \frac{1}{t^{N-1}}%
\int_{0}^{t}s^{N-1}a_{i}\left( s\right) f_{i}\left( u_{1}^{0}\left( s\right)
,...,u_{m}^{0}\left( s\right) \right) ds\right) ^{1/\left( p-1\right) }dt \\
&\leq &z\left( 0\right) +\int_{0}^{r}\left( \frac{1}{t^{N-1}}%
\int_{0}^{t}s^{N-1}\overset{m}{\underset{i=1}{\Sigma }}a_{i}\left( s\right) 
\overset{m}{\underset{i=1}{\Sigma }}f_{i}\left( z\left( s\right)
,...,z\left( s\right) \right) ds\right) ^{1/\left( p-1\right) }dt=z\left(
r\right) .
\end{eqnarray*}%
Thus $u_{i}^{1}\left( r\right) \leq z\left( r\right) $ ($i=1,...,m$).
Similar arguments show that%
\begin{equation*}
u_{i}^{k}\left( r\right) \leq z\left( r\right) \text{ for all }r\in \left[
0,\infty \right) \text{ and }k\geq 1\text{.}
\end{equation*}%
Then, we may assume 
\begin{equation*}
\left( u_{1}\left( \left\vert x\right\vert \right) ,...,u_{m}\left(
\left\vert x\right\vert \right) \right) :=\left( \underset{k\rightarrow
\infty }{\lim }u_{1}^{k}\left( \left\vert x\right\vert \right) ,...,\underset%
{k\rightarrow \infty }{\lim }u_{m}^{k}\left( \left\vert x\right\vert \right)
\right) \text{, for every }x\in \mathbb{R}^{N}
\end{equation*}%
is an entire radial solution of system (\ref{11}).

Now, let $\left( u_{1},...,u_{m}\right) $ be any non-negative non-trivial
entire radial solution of (\ref{11}) and suppose that $a_{i}$ ($i=1,...,m$)
satisfies (\ref{12}). Since $u_{i}$ ($i=1,...,m$) is nontrivial and
non-negative, there exists $R>0$ so that $u_{i}\left( R\right) >0$. Since $%
u_{i}^{\prime }\geq 0$, we get $u_{i}\left( r\right) \geq u_{i}\left(
R\right) $ for $r\geq R$ and thus from%
\begin{equation*}
u_{i}\left( r\right) =u_{i}\left( 0\right) +\int_{0}^{r}\frac{1}{t^{N-1}}%
\left( \int_{0}^{t}s^{N-1}a_{i}\left( s\right) f_{i}\left( u_{1}\left(
s\right) ,...,u_{m}\left( s\right) \right) ds\right) ^{1/\left( p-1\right)
}dt,
\end{equation*}%
we obtain

\begin{equation*}
\left\{ 
\begin{array}{l}
u_{i}\left( r\right) =u_{i}\left( 0\right) +\int_{0}^{r}\left( \frac{1}{%
t^{N-1}}\int_{0}^{t}s^{N-1}a_{i}\left( s\right) f_{i}\left( u_{1}\left(
s\right) ,...,u_{m}\left( s\right) \right) ds\right) ^{1/\left( p-1\right)
}dt\text{ } \\ 
\text{ \ \ \ \ \ \ }\geq u_{i}\left( R\right) +f_{i}^{1/\left( p-1\right)
}\left( u_{1}\left( R\right) ,...,u_{m}\left( R\right) \right)
\int_{R}^{r}\left( \frac{1}{t^{N-1}}\int_{R}^{t}s^{N-1}a_{i}\left( s\right)
ds\right) ^{1/\left( p-1\right) }dt\rightarrow \infty \text{ as }%
r\rightarrow \infty , \\ 
\text{for all }i=1,...,m.%
\end{array}%
\right. 
\end{equation*}%
Conversely, if $f_{i}$ ($i=1,...,m$) satisfy (C1)-(C3) and $\left(
w_{1},..,w_{m}\right) $ is a nonnegative entire large solution of (\ref{11}%
), then $w_{i}$ satisfy%
\begin{eqnarray*}
\left( p-1\right) w_{1}^{\prime }\left( r\right) ^{p-2}w_{1}^{\prime \prime
}+\frac{N-1}{r}w_{1}^{\prime }\left( r\right) ^{p-1} &=&a_{1}\left( r\right)
f_{1}\left( w_{1},...,w_{m}\right) , \\
&&... \\
\left( p-1\right) w_{m}^{\prime }\left( r\right) ^{p-2}w_{m}^{\prime \prime
}+\frac{N-1}{r}w_{m}^{\prime }\left( r\right) ^{p-1} &=&a_{m}\left( r\right)
f_{i}\left( w_{1},...,w_{m}\right) .
\end{eqnarray*}%
Then, using the monotonicity of \textit{\ }$r^{\frac{p\left( N-1\right) }{p-1%
}}\underset{j=1}{\overset{m}{\Sigma }}a_{j}\left( r\right) $\textit{\ } we
can apply similar arguments used in obtaining \textbf{Theorem \ref{1} }to
get 
\begin{eqnarray*}
\left( \overset{m}{\underset{i=1}{\Sigma }}w_{i}\left( r\right) \right)
^{\prime } &\leq &\left[ Cr^{\frac{p\left( 1-N\right) }{p-1}}+\frac{pm^{p-1}%
}{p-1}\overset{m}{\underset{i=1}{\Sigma }}a_{i}\left( r\right) F\left( 
\overset{m}{\underset{i=1}{\Sigma }}w_{i}\left( r\right) \right) \right]
^{1/p} \\
&\leq &\sqrt[p]{C}r^{\frac{1-N}{p-1}}+\sqrt[p]{\frac{pm^{p-1}}{p-1}\overset{m%
}{\underset{i=1}{\Sigma }}a_{i}\left( r\right) }\left[ F\left( \overset{m}{%
\underset{i=1}{\Sigma }}w_{i}\left( r\right) \right) \right] ^{1/p},
\end{eqnarray*}%
and hence as in (\ref{111}), we get

\begin{eqnarray*}
&&\int_{\overset{m}{\underset{i=1}{\Sigma }}w_{i}\left( R\right) }^{\overset{%
m}{\underset{i=1}{\Sigma }}w_{i}\left( r\right) }\left[ F\left( t\right) %
\right] ^{-1/p}dt \\
&\leq &\frac{\sqrt[p]{C}}{\left[ F\left( \overset{m}{\underset{i=1}{\Sigma }}%
w_{i}\left( R\right) \right) \right] ^{1/p}}\int_{R}^{r}t^{\frac{1-N}{p-1}%
}dt+\sqrt[p]{\frac{1}{2}\frac{pm^{p-1}}{p-1}}\left(
\int_{R}^{r}t^{1+\varepsilon }\left( \overset{m}{\underset{i=1}{\Sigma }}%
a_{i}\left( t\right) \right) ^{2/p}dt+\int_{R}^{r}t^{-1-\varepsilon
}dt\right) \\
&\leq &\sqrt[p]{C}\left[ F\left( \overset{m}{\underset{i=1}{\Sigma }}%
w_{i}\left( R\right) \right) \right] ^{-1/p}\frac{\left( p-1\right) R^{\frac{%
p-N}{p-1}}}{p-N}+\sqrt[p]{\frac{1}{2}\frac{pm^{p-1}}{p-1}}\left(
\int_{R}^{r}t^{1+\varepsilon }\left( \overset{m}{\underset{i=1}{\Sigma }}%
a_{i}\left( t\right) \right) ^{2/p}dt+\frac{1}{\varepsilon R^{\varepsilon }}%
\right) \\
&\leq &C_{R}+\int_{R}^{r}t^{1+\varepsilon }\left( \overset{m}{\underset{i=1}{%
\Sigma }}a_{i}\left( t\right) \right) ^{2/p}dt\text{,}
\end{eqnarray*}%
where 
\begin{equation*}
C_{R}=\sqrt[p]{C}\left[ F\left( \overset{m}{\underset{i=1}{\Sigma }}%
w_{i}\left( R\right) \right) \right] ^{-1/p}\frac{\left( p-1\right) R^{\frac{%
p-N}{p-1}}}{p-N}+\frac{1}{\varepsilon R^{\varepsilon }}.
\end{equation*}%
Passing to the limit as $r\rightarrow \infty $, we find that $a_{j}$ ($%
j=1,...,m$) satisfies (\ref{13}). \rule{5pt}{5pt}

\begin{remark}
\textit{If (C1)-(C3) are satisfied then}%
\begin{equation*}
\int_{1}^{\infty }[\int_{0}^{s}f_{i}\left( t,...,t\right)
dt]^{-1/p}ds=\infty \text{, }i=1,...,m\text{.}
\end{equation*}
\end{remark}

\begin{remark}
(see \cite{CD}) If (C1)-(C2) and 
\begin{equation*}
\int_{1}^{\infty }\left( \underset{i=1}{\overset{m}{\Sigma }}f_{i}\left(
s,...,s\right) \right) ^{-1/\left( p-1\right) }ds=\infty ,
\end{equation*}%
are satisfied, then%
\begin{equation*}
\int_{1}^{\infty }\left( \int_{0}^{t}\underset{i=1}{\overset{m}{\Sigma }}%
f_{i}\left( s,...,s\right) \right) ^{-1/p}dsdt=\infty .
\end{equation*}
\end{remark}

\end{document}